\documentclass[a4paper,12pt]{article}
\usepackage[english]{babel}
\usepackage{amsmath}
\usepackage{amssymb}
\usepackage{tabls}
\usepackage[pdftex]{graphicx}
\usepackage{calc}
\usepackage{xcolor}
\usepackage[T1]{fontenc}
\usepackage[latin1]{inputenc}
\title{Dimensions of Crystalline Graded Rings}
\author{Tim Neijens\\University of Antwerp\\ \texttt{tim.neijens@gmail.com} \and Freddy Van Oystaeyen \\ University of Antwerp \\ \texttt{fred.vanoystaeyen@ua.ac.be}}

\newcommand{\blok}{\hfill \Box}
\newcommand{\CGR}{\mathop \diamondsuit \limits_{\sigma ,\alpha}}
\newtheorem{defi}{Definition}[section]
\newtheorem{opm}[defi]{Remark}
\newtheorem{prop}[defi]{Proposition}
\newtheorem{gev}[defi]{Corollary}
\newtheorem{st}[defi]{Theorem}
\newtheorem{lem}[defi]{Lemma}

\begin{document}
\maketitle

\begin{abstract}
The global dimension of a ring governs many useful abilities.  For example, it is semi-simple if the global dimension is $0$, hereditary if it is $1$ and so on.  We will calculate the global dimension of a Crystalline Graded Ring, as defined in the paper by E. Nauwelaerts and F. Van Oystaeyen, \cite{NVO6}.  We will apply this to derive a condition for the Crystalline Graded Ring to be semiprime.  In the last section, we give a little bit of attention to the Krull-dimension.
\end{abstract}

\section{Preliminaries}

{\defi \label{def1}\textbf{Pre-Crystalline Graded Ring}\\
Let $A$ be an associative ring with unit $1_A$.  Let $G$ be an arbitrary group.  Consider an injection $u: G \rightarrow A$ with $u_e = 1_A$, where $e$ is the neutral element of $G$ and $u_g \neq 0$,  $\forall g \in G$.  Let $R \subset A$ be an associative ring with $1_{R}=1_A$.  We consider the following properties:
\begin{description} 
	\item[(C1)]\label{def2} $A = \bigoplus_{g \in G} R u_g$.
	\item[(C2)]\label{def3} $\forall g \in G$, $R u_g = u_g R$ and this is a free left $R$-module of rank $1$.
	\item[(C3)]\label{def4} The direct sum $A = \bigoplus_{g \in G} R u_g$ turns $A$ into a $G$-graded ring with $R = A_e$.
\end{description}
We call a ring $A$ fulfilling these properties a \textbf{Pre-Crystalline Graded Ring}.}\\

{\prop \label{def5} With conventions and notation as in Definition \ref{def1}:
\begin{enumerate}
	\item For every $g \in G$, there is a set map $\sigma_g : R \rightarrow R$   defined by: $u_g r = \sigma_g(r)u_g$ for $r \in R$.  The map $\sigma_g$ is in fact a surjective ring morphism.  Moreover, $\sigma_e = \textup{Id}_{R}$.
	\item There is a set map $\alpha : G \times G \rightarrow R$ defined by $u_g u_h = \alpha(g,h)u_{gh}$ for $g,h \in G$.  For any triple $g,h,t \in G$ the following equalities hold:
		\begin{eqnarray}
		\alpha(g,h)\alpha(gh,t)&=&\sigma_g(\alpha(h,t))\alpha(g,ht) \label{def6},\\
		\sigma_g(\sigma_h(r))\alpha(g,h)&=& \alpha(g,h)\sigma_{gh}(r) \label{def7}.
		\end{eqnarray}
	\item $\forall g \in G$ we have the equalities $\alpha(g,e) = \alpha(e,g) = 1$ and $\alpha(g,g^{-1}) = \sigma_g(\alpha(g^{-1},g)).$
\end{enumerate}
}
\begin{flushleft}\textbf{Proof}\end{flushleft} See \cite{NVO6}. $\blok$\\

{\prop Notation as above, the following are equivalent:
\begin{enumerate}
	\item $R$ is $S(G)$-torsionfree.
	\item $A$ is $S(G)$-torsionfree.
	\item $\alpha(g,g^{-1})r=0$ for some $g \in G$ implies $r = 0$.
	\item $\alpha(g,h)r=0$ for some $g,h \in G$ implies $r = 0$.
	\item $R u_g = u_g R$ is also free as a right $R$-module with basis $u_g$ for every $g \in G$.
	\item for every $g \in G$, $\sigma_g$ is bijective hence a ring automorphism of $R$.
\end{enumerate}
}
\begin{flushleft}\textbf{Proof}\end{flushleft} See \cite{NVO6}. $\blok$\\

{\defi Any $G$-graded ring $A$ with properties \textbf{(C1),(C2),(C3)}, and which is $G(S)$-torsionfree is called a \textbf{crystalline graded ring}.  In case $\alpha(g,h) \in Z(R)$, or equivalently $\sigma_{gh}=\sigma_g \sigma_h$, for all $g,h \in G$, then we say that $A$ is \textbf{centrally crystalline}.}\\

{\lem \label{def9}Let $R \CGR G$ be a pre-crystalline graded ring, $x \in R$, $g,h \in G$.  $R$ is a domain, and define $K$ to be the quotient field of $R$.  Then
\begin{enumerate}
  \item $u_g^{-1} = u_{g^{-1}}\alpha^{-1}(x,x^{-1})=\alpha^{-1}(x^{-1},x)u_{x^{-1}}$.
	\item $\sigma_g^{-1}(x)u_g^{-1} = u_g^{-1}x$.
	\item $\sigma_{hg}^{-1}[\alpha(h,g)]=\sigma_g^{-1}[\sigma_h^{-1}(\alpha(h,g))]$.
	\item $\sigma_g^{-1}[\alpha(g,g^{-1}h)]=\alpha^{-1}(g^{-1}, h)\sigma_g^{-1}[\alpha(g, g^{-1})]$.
\end{enumerate}
}
\begin{flushleft}\textbf{Proof}\end{flushleft}(inverses are defined in $K$ or $K \CGR G$)
\begin{enumerate}
	\item Just calculate the product and use that in an associative ring the left and right inverse coincide.
	\item Let $g,h \in G, x \in A$:
	\begin{align*}
	&\sigma_g[\sigma_h(x)]\alpha(g,h)=\alpha(g,h)\sigma_{gh}(x)\\
	\Rightarrow &\sigma_g[\sigma_{g^{-1}}(x)]\alpha(g,g^{-1})=\alpha(g,g^{-1})x\\
	\Rightarrow &\sigma_{g^{-1}}(x)\sigma_g^{-1}(\alpha(g,g^{-1}))=\sigma_{g^{-1}}(\alpha(g,g^{-1}))\sigma_{g}^{-1}(x)\\
	\Rightarrow &\sigma_g^{-1}(x) = \sigma_g^{-1}[\alpha^{-1}(g,g^{-1})]\sigma_{g^{-1}}(x)\sigma_g^{-1}[\alpha(g,g^{-1})].
	\end{align*}
	So
	\begin{align*}
	\sigma_g^{-1}(x)u_g^{-1} &= \sigma_g^{-1}[\alpha^{-1}(g,g^{-1})]\sigma_{g^{-1}}(x)\sigma_g^{-1}[\alpha(g,g^{-1})]\alpha^{-1}(g^{-1},g)u_{g^{-1}}\\
	\ &= \sigma_g^{-1}[\alpha^{-1}(g,g^{-1})]\sigma_{g^{-1}}(x)\alpha(g^{-1},g)\alpha^{-1}(g^{-1},g)u_{g^{-1}}\\
	\ &= \sigma_g^{-1}[\alpha^{-1}(g,g^{-1})]\sigma_{g^{-1}}(x)u_{g^{-1}}\\
	\ &= \sigma_g^{-1}[\alpha^{-1}(g,g^{-1})]u_{g^{-1}}x\\
	\ &= \alpha^{-1}(g,g^{-1})u_{g^{-1}}x\\
	\ &= u_g^{-1}x.
	\end{align*}
	\item Let $g,h \in G, x \in A$:
	\begin{align*}
	&\sigma_h[\sigma_g(x)]\alpha(h,g) = \alpha(h,g)\sigma_{hg}(x)\\
	\Rightarrow &\sigma_h[\sigma_g(\sigma_{hg}^{-1}(\alpha(h,g)))]\alpha(h,g) = \alpha(h,g)\sigma_{hg}(\sigma_{hg}^{-1}(\alpha(h,g)))\\
	\Rightarrow &\sigma_{hg}^{-1}[\alpha(h,g)] = \sigma_g^{-1}[\sigma_h^{-1}(\alpha(h,g))].
	\end{align*}
	\item Let $g,h \in G$:
	\[\alpha(g,g^{-1})\alpha(e,h)=\sigma_g[\alpha(g^{-1},h)]\alpha(g,g^{-1}h).\]
	\ $\blok$
\end{enumerate}

\section{Global Dimension}

{\st \label{dim1} Let $R,S$ be rings with $R \subseteq S$ such that $R$ is an $R$-bimodule direct summand of $S$, then $\textup{r gld}R \leq \textup{r gld}S + \textup{pd} S_R$.}\\
\textbf{Proof} See \cite{MR},p. 237. $\blok$\\

{\st \label{dim2} Let $R$ be a ring, $G$ a finite group with $|G|$ a unit in $R$ and $A = R \CGR G$ a pre-crystalline graded ring with $u_g$ units.  Let $M$ be any right $A$-module.  Then:
\begin{enumerate}
	\item If $N \triangleleft M_A$ and $N$ is a direct summand of $M$ as an $R$-module, then $N$ is a direct summand over $A$.
	\item $\textup{pd}M_{R} = \textup{pd}M_A$.
	\item $\textup{r gld}R = \textup{r gld}A$.
\end{enumerate}
}
\begin{flushleft}\textbf{Proof}\end{flushleft}
\begin{enumerate}
	\item \label{dim3} Let $\pi : M \rightarrow N$ be the $R$-module splitting morphism.  Define the map $\lambda$ by
	\[\lambda: M \rightarrow N :m \mapsto |G|^{-1}\sum_{g \in G}\pi(m u_g) u_g^{-1}.\]
	\textbf{$\lambda$ is well-defined} : trivial.\\
	\textbf{$\lambda$ is the identity on $N$} : let $k \in N$:
	\begin{align*}
	\lambda(k) &=|G|^{-1}\sum_{g \in G} \pi(k u_g)u_g^{-1}\\
	&= |G|^{-1}\sum_{g\in G}k = k.
	\end{align*}
	\textbf{$\lambda$ is $A$-linear} : Let $m \in M, a\in A$:
	\begin{align*}
	\lambda(ma) =& |G|^{-1} \sum_{g\in G} \pi(mau_g)u_g^{-1}\\
	=& |G|^{-1} \sum_{g\in G} \pi\left[m\left(\sum_{h \in G} t_h u_h\right)u_g\right]u_g^{-1}\\
	=& |G|^{-1} \sum_{g,h\in G} \pi\left(m t_h u_h u_g\right)u_g^{-1}\\
	{}^{(\textup{Lemma }\ref{def9}(2))}=& |G|^{-1} \sum_{g,h\in G} \pi\left(m u_h u_g\right)u_g^{-1}\sigma_h^{-1}(t_h)\\
		\end{align*}
	\begin{align*}
	=& |G|^{-1} \sum_{g,h\in G} \pi\left(m \alpha(h,g)u_{hg}\right)u_g^{-1}\sigma_h^{-1}(t_h)\\
	=& |G|^{-1} \sum_{g,h\in G} \pi\left(m u_{hg}\right)\sigma_{hg}^{-1}[\alpha(h,g)]u_g^{-1}\sigma_h^{-1}(t_h)\\
	{}^{(\textup{Lemma }\ref{def9}(3))} =&|G|^{-1} \sum_{g,h\in G} \pi\left(m u_{hg}\right)\sigma_g^{-1}[\sigma_h^{-1}(\alpha(h,g))]u_g^{-1}\sigma_h^{-1}(t_h)\\
	{}^{(\textup{Lemma }\ref{def9}(2))} =&|G|^{-1} \sum_{g,h\in G} \pi\left(m u_{hg}\right)u_g^{-1}\sigma_h^{-1}[\alpha(h,g)]\sigma_h^{-1}(t_h)\\
	{}^{(x=hg)}=& |G|^{-1} \sum_{h\in G} \sum_{x \in G}\pi\left(m u_x\right)u_{h^{-1}x}^{-1}\sigma_h^{-1}[\alpha(h,h^{-1}x)]\sigma_h^{-1}(t_h)\\
	{}^{(\textup{Lemma \ref{def9}(4)})}=& |G|^{-1} \sum_{h\in G} \sum_{x \in G}\pi\left(m u_x\right)[\alpha^{-1}(h^{-1},x)u_{h^{-1}}u_x]^{-1}\cdot\\
	&\qquad \qquad \qquad \alpha^{-1}(h^{-1},x)\sigma_h^{-1}[\alpha(h,h^{-1})]\sigma_h^{-1}(t_h)\\
	=& |G|^{-1} \sum_{h\in G} \sum_{x \in G}\pi\left(m u_x\right)u_x^{-1} u_{h^{-1}}^{-1}\sigma_h^{-1}[\alpha(h,h^{-1})]\sigma_h^{-1}(t_h)\\
	=& |G|^{-1}\sum_{h \in G}\sum_{x \in G} \pi(mu_x)u_x^{-1}u_h \sigma_h^{-1}(t_h)\\
	=& |G|^{-1}\sum_{x \in G}\pi(m u_x)u_x^{-1}\sum_{h \in G}t_h u_h\\
	=& \lambda(m)\cdot a.
	\end{align*}
\item\label{dim4}
Suppose $M_{R}$ is projective and
\[0 \rightarrow N \rightarrow F \rightarrow M \rightarrow 0\]
is a short exact sequence of $A$-modules with $F$ free, then the sequence splits over $R$ and hence over $A$ by (\ref{dim3}).  So $M_A$ is also projective.  Furthermore, $A_{R}$ is free.  It now follows that an $A$-projective resolution of any module $M_A$ is also an $R$-projective resolution that terminates when a kernel is, equally, $R$-projective or $A$-projective, so $\textup{pd}M_{R} = \textup{pd}M_A$.
\item
Any $A$-module is naturally an $R$-module.  So, since $\textup{pd}M_{R} = \textup{pd}M_A$, we find
\begin{eqnarray*}
\textup{r gld}A &=& \textup{sup}\left\{\textup{pd}M_A | M_A \textup{ right } A-\textup{module}\right\}\\
\ &\leq& \textup{sup}\left\{\textup{pd}M_{R} | M_{R} \textup{ right } R-\textup{module}\right\}\\
\ &=& \textup{r gld}R.
\end{eqnarray*}

So by Theorem \ref{dim1}:
\begin{eqnarray*}
\textup{r gld}R &\leq& \textup{r gld}A + \textup{pd}A_{R}\\
\ &\mathop  = \limits^{(\ref{dim4})}& \textup{r gld}A + \textup{pd}A_A\\
\ &=& \textup{r gld}A.
\end{eqnarray*}
And in conclusion $\textup{r gld}R = \textup{r gld}A$.$\blok$
\end{enumerate}

The following result is well-known:

{\lem \label{dim5} Let $S$ be an Ore set for $R$ and suppose there is no $S$-torsion.  Let $\{s_1, \ldots, s_n\} \subset S$, then $\exists s \in S \cap \bigcap_{i =1}^{n} R s_i$.}\\
\textbf{Proof} By induction.  Let us take $s_1$, $\exists t_1 \in S^{-1}R$ such that $t_1 s_1 = 1$.  Then of course we find $q_1 \in S$ such that $q_1t_1 \in R$.  This means that $q_1 = st_1 s_1 \in R s_1$, and $q_1 \in S$.  Now we try to do the same for the other $s_i$.  We apply the left Ore condition on $q_1 \in S \subset R$ and $s_2 \in S$.  We now find $v_2 \in R$ and $q_2 \in S$ such that $v_2 s_2 =  q_2 q_1$. $\blok$\\

{\lem Let $A = R \CGR G$ be crystalline graded, then the set of regular elements in $R$, $\textup{reg}R$, is a subset of $\textup{reg}A$, the regular elements of $A$.  Furthermore, if $R$ is semiprime Goldie, $\textup{reg}R$ is a left (and right) Ore set in $A$.  We have
\[\left(\textup{reg}R\right)^{-1}A = \bigoplus_{g \in G}Q_{\textup{cl}}(R)u_g.\]}
\begin{flushleft}\textbf{Proof}\end{flushleft} For the first part, take $a \in \textup{reg}R$, $x = \sum_{g \in G} x_g u_g$ and suppose $ax=0$, then $\sum_{g \in G} ax_gu_g =0$.  This implies $a x_g=0$ $\forall g \in G$, and this means $x_g, \forall g \in G$.  Suppose $xa = 0$, then $\sum_{g \in G}x_g u_g a = 0$.  This implies $x_g \sigma_g(a)u_g = 0$, or $x_g \sigma_g(a) = 0, \forall g \in G$.  Since $\textup{reg}R$ is invariant under $\sigma_g, \forall g \in G$, we again find $x_g = 0, \forall g \in G$.  So we have proven $\textup{reg}R \subset \textup{reg}A$.\\
By Goldie's Theorem, we know that $\textup{reg}R$ is an Ore set in $R$.  We first need to prove that $S=\textup{reg}R$ satisfies the left Ore condtion for $A$.  We need that $\forall r \in R$, $s \in S$ we can find $r' \in R, s' \in S$ such that $s'r = r's$.  Let $r = \sum_{g \in G} a_g u_g$.  Since $S$ is left Ore for $R$, we can find $\forall g \in G$ elements $a'_g \in R$ and $s_g \in S$ such that $a'_g \sigma_g(s)=s_g a_g$.  Now, we find $s' \in S \cap \bigcap_{g \in G} R s_g$ from Lemma \ref{dim5}, in other words, we find $s' \in S$ and $v_g \in R$ such that $\forall g \in G$ $s' = v_g s_g$.  Now set $\forall g \in G$, $b_g = v_g a'_g$, and set $r' = \sum_{g \in G}b_g u_g$.  Then $r' s = s' r$.  The right Ore condition is similar.  The third assertion is now clear. $\blok$

{\st \label{dim6} Let $A$ be crystalline graded over $R$, $R$ a semiprime Goldie ring.  Assume $\textup{char}R$ does not divide $|G|$, then $A$ is semiprime Goldie.}\\
\textbf{Proof} Since $A$ is crystalline graded, the elements $\alpha(g,h), g,h \in G$ are regular elements.  Denote $S = \textup{reg}R$.  Since $R$ is semiprime Goldie, $S^{-1}R$ is semisimple Artinian.  This implies that from Theorem \ref{dim2}, $S^{-1}A$ is semisimple Artinian, in particular, it is Noetherian.  Let $I$ be an ideal in $A$, and consider $(S^{-1}A) I$.  Claim: this is an ideal.  Let $s \in S$ and consider the following chain:
\[(S^{-1}A) I\subset (S^{-1}A) I s^{-1} \subset (S^{-1}A) I s^{-2} \subset \ldots.\]
This implies that $(S^{-1}A) Is^{-n} = (S^{-1}A) Is^{-m}$, $m >n$, and so $(S^{-1}A) I = (S^{-1}A) I s^{n-m}$, and so we find $(S^{-1}A) I(S^{-1}A)\subset (S^{-1}A) I$, or $(S^{-1}A) I$ is an ideal in $S^{-1}A$.  If $J$ is the nilradical of $A$ then $(S^{-1}\cdot J)^n = S^{-1} \cdot J^n$ follows.  For some $n$ we have that $(S^{-1}\cdot J)^n = 0$ in the semisimple Artinian ring $S^{-1}A$, thus $S^{-1}A \cdot J = 0$ and $J=0$. $\blok$\\

{\gev If $A$ is crystalline graded with $D$ a Dedekind domain, $\textup{char}D$ does not divide $|G|$, then $A$ is semiprime.}\\

{\prop In the situation of Theorem \ref{dim6}, prime ideals of $S^{-1}A$ intersect in prime ideals of $A$, where $S = \textup{reg}R$.}\\
\textbf{Proof} Let $P$ be a prime of $S^{-1}A$, then $P \cap Q$ is an ideal such that for $IJ \subset P \cap A$, $I$ and $J$ ideals of $A$, we have $S^{-1}A\cdot IJ \subset P$ hence $(S^{-1}A\cdot I)(S^{-1}A\cdot J) \subset P$, or $S^{-1}A\cdot I \in P$ if $S^{-1}A \cdot J \not\subset P$.  Thus $I \subset P\cap A$ if $J \not\subset P\cap A$ and conversely. $\blok$\\

{\opm The situation of Theorem \ref{dim6} arises when $A$ is centrally crystalline graded over the semiprime Goldie ring $R$ with $\textup{char}R$ does not divide $|G|$, such that $A$ (or $R$) is a P.I. ring.}

\section{Krull Dimension}

{\prop Let $A$ be crystalline graded over $D$, $D$ a Dedekind domain.  Then the (Krull-)dimension of $A$ is smaller than or equal to $2$.}\\
\textbf{Proof} Consider the set $F = \{I \triangleleft A | I \cap D =0\}$ ordered by inclusion.  If it is nonempty, then there is a maximal element for this family, say $P$.  Suppose $IJ \subset P$, with $P \not\subset P+I$, $P \not\subset P+J$.  Then $0\neq d_1 \in P + I \cap D$ and $0 \neq d_2 \in P+J \cap D$.  This implies $0 \neq d_1d_2 \in P$, contradiction.  So if $F \neq \emptyset$, there always exists a prime ideal $P$ in $A$ with $P \cap D = 0$.\\
Denote $S = D \backslash\{0\}$. Suppose that $0 \neq Q \subset P$, $Q$ a prime ideal in $A$.  Then, since $S^{-1}A$ is Artinian semisimple (Theorem \ref{dim2}), we find that $S^{-1}Q = S^{-1}P$ since they are both primes ($Q \cap D \neq 0 \neq P \cap D$).  Now let $y \in P \backslash Q$.  Then $y \in S^{-1}P = S^{-1}Q$.  This means $\exists d \in S$ such that $dy \in Q$.  So if we set $d'= \prod_{g \in G} \sigma_g(d)$ then $d'y\in Q$.  Since $d' \in Z(A)$ we find $d'Ay \subset Q$ and since $y \notin Q$ we see that $d' \in Q$ or $Q \cap D \neq 0$.  Contradiction.  We have established that two prime ideals that don't intersect $D$ cannot contain each other.\\
Suppose there exists a prime ideal $M$ of $A$ with $M \cap D \neq 0$.  This means $A/M$ is Artinian, and prime, in other words it is a simple ring, or $M$ is a maximal ideal.  We find that a maximal chain of prime ideals always is of the form
\[0 \subset P \subset M \subset A,\]
where $P \cap D = 0$ and $Q \cap D \neq 0$. $\blok$\\

\end{document}